\documentclass[10pt]{article}

\usepackage{graphicx}
\usepackage{graphicx,psfrag}
\usepackage{amsmath,amssymb,amsthm}
 
\def \vs {\vskip 0.2cm}

\def \n {\noindent}

\def \bN {{\mathbb  N}}

\def \bE {{\mathbb E}}

\def \bC {{\mathbb C}}

\def \CD {{\cal D}}

\def \CT {{\cal T}}

\def   \l  {{  \lambda}}

\def \CP {{\cal P}}

\def \Cum{{\rm Cum}}

\def \s {\sigma}

\begin{document}
\title{Enumeration of tree-type diagrams assembled from 
oriented chains of edges\footnote{{\bf MSC: 05A15, 05C30, 60B20} 
}
}

\author{O. Khorunzhiy\\ Universit\'e de Versailles - Saint-Quentin \\45, Avenue des Etats-Unis, 78035 Versailles, FRANCE\\
{\it e-mail:} oleksiy.khorunzhiy@uvsq.fr}
\maketitle
\begin{abstract}

We study a family of tree-type diagrams that arise in studies of 
the cumulant expansion in discrete Erd\H os-R\'enyi random matrix models. 
Using 
a version of the Pr\" ufer code, 
we  obtain  an explicit expression for the number of  
 tree-type diagrams assembled from
$k$   oriented chains of $q$ edges.
Using this modified Pr\" ufer codification, we get 
an explicit expression  for  sum overs weighted 
tree-type diagrams with a 
weight depending on multiplicity of edges. 
We describe similar results for 
 tree-type diagrams assembled from  chains that are not necessarily regular.

  \end{abstract}


{\it Key words:}  Pr\" ufer code, weighted trees

\section{Random graphs, walks  and  tree-type diagrams}

Let us consider a family of jointly independent identically distributed Bernoulli random variables  
$$
a_{ij} = 
\begin{cases}
 1, & \text{with probability   $p_n$} , \\
0,  & 
\text {with probability   $1-p_n$,} 
\end{cases}  \quad 1\le i< j\le n.
\eqno (1) 
$$
A real symmetric $n$-dimensional matrix $A$ 
with elements 
$
(A)_{ij}= 
 a_{ij}$, $1\le i<j\le n$ and $(A)_{ii}=0$
 can be regarded as  adjacency matrix of random graphs of the 
  Erd\H os-R\'enyi ensemble \cite{B}.
In this representation, the sum 
$$
Y_n^{(q)} = \sum_{i,j = 1}^n \left( A^q\right)_{ij}
=\sum_{i_1, i_2, \dots, i_{q+1} =1}^n a_{i_1i_2} a_{i_2i_3} \dots a_{i_qi_{q+1}}
\eqno(2) 
$$
can be regarded as 
 the  number  of 
$q$-step walks
on a realization of the  Erd\H os-R\'enyi random graph. 
In paper \cite{K-08}, the following 
 cumulant expansion  has been considered
$$
f^{(n,\beta)}(g) = \log \bE \exp\left\{ g 
Y_n^{(q)}\right\}
= \sum_{k\ge 1} {g^k\over k!} \Cum_k(Y^{(q)}_n),
\eqno (3)
$$
where 
$\bE$ denotes the mathematical expectation with respect to the measure generated by
random variables (1) and $\Cum_k(Y_n^{(q)})$ denotes the $k$-th cumulant of $Y_n^{(q)}$.  It can be considered also as a definition of cumulants of $Y_n^{(q)}$ \cite{L}. 
The last equality of (3) is known in mathematical physics as the cumulant expansion of the Helmholtz free energy $f^{(n,\beta)}(g)$. Variables $a_{ij}$ (1) are bounded and therefore the series of (3) is absolutely convergent for any finite $n$ and sufficiently small $g$. 

\begin{figure}[htbp]
\centerline{\includegraphics
[height=6cm, width=16cm]{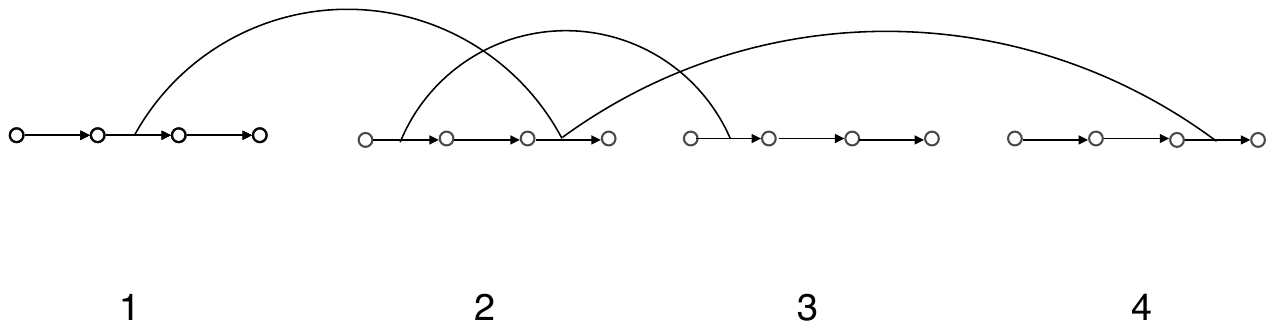}}
\caption{{ A diagram $D_k^{(q)} $ with $q=3$ and $k=4$}}
\end{figure}

\noindent In paper \cite{K-08}, it  is shown that 
in the limiting transition when $p= c/n$ and  $n, c\to\infty$,
$$
{1\over n c} \Cum_k
\left( {1\over c^q}Y^{(q)}_n\right)
= 2^{k-1}\ d^{(q)}_k (1+ o(1)), \quad n, c\to\infty,
\eqno (4)
$$
where 
 $d^{(q)}_k$ is given by  the number of all possible connected tree-type diagrams 
assembled from $k$  linear  graphs with $q$ oriented edges.

Tree-type diagrams 
can be determined as follows.  
Let us consider
$k$ elements 
$\l_q^{(i)}$, \mbox{$1\le i\le k$,} each element is represented
by a linear graph of $q$ ordered oriented edges. We refer to elements 
$\l_q^{(i)}$ as to $q$-regular chains. 
These elements are connected by $k-1$ arcs. Each arc  joins   an edge of one element with that of to another one. 
We say that these edges represent the feet of the arc. The arcs are drawn in such way that each element $\l_q$ has at least one edge that serves as a foot for one or several arcs. 
We say that the ensemble of $\l_q$-elements together with arcs represents
a diagram $D_k^{(q)}$. 
The diagrams obtained are connected. On figure 1, we give an example 
of $D_k^{(q)}$ with $q=3$ and $k=4$. 
It is clear that identifying  edges joined by an arc and forgetting the orientation of edges of $\l_q$,
we get a tree graph. Therefore  we say that $D_k^{(q)}$ is a  tree-type diagram
and denote by $\CD_k^{(q)}$ the ensemble of  all possible tree-type diagrams.

In paper \cite{K-08}, it is shown that a kind of exponential generating function $H$ of numbers $d_k^{(q)}$, $ k\ge 0$ verifies a differential equation that can be solved to get a kind of P\'olya equation for the power of $H$. Basing on this result, an explicit expression has been obtained for 
$d_k^{(q)}$ in the particular case when $q=2$,
$$
d^{(2)}_k = 2^k (k+1)^{k-2}, \quad k\in \bN.
\eqno (5)
$$
The number $d_k^{(2)}$ is known as 
the number of directed trees with $k+1$ vertices and  $k$ 
distinct edge labels \cite{BBR,OEIS} 
that can be determined as graphs dual  to tree-type diagrams $D_k^{(2)}$ 
(see \cite{K-08} for more details). The sequence 
$\{d_k^{(2)}, k\ge 1\}
= \{1,4, 32, 400, \dots\}$ is referred in OEIS as the discriminants of Chebyshev polynomials \cite{OEIS}.

In the present paper we consider the  general case of $q\ge 2$.  
The  main result of this paper is as follows.  

{\bf Theorem 1.} {\it Given $q\ge 2$, the number of all possible tree-type diagrams $d_k^{(q)}= | \CD_k^{(q)}|$ 
is given by expression
$$
d^{(q)}_k = q^k \big((q-1)k+1\big)^{k-2}, \quad k \in \bN.
\eqno (6) 
$$
Remark.} Theorem 1 is also valid in the trivial case of $q=1$. 
We did not find any reference for the sequence (6) in the literature. 
Also we did not find any reference in OEIS  for  (6) with given $q$, 
say $q=3$. 
 
 \vs 
To prove Theorem 1, we 
develop an enumeration technique based on a version of the Pr\"ufer codification procedure.
As a further use of our technique, 
we consider two generalizations of the main result (6). 
First, 
we  get an explicit expression  for the total sum of  weighted $q$-regular tree-type diagrams such that the weight of diagram $W(D_k^{(q)})$ depends on the number of multiple edges in $\CD_k^{(q)}$. 
The second generalization is given by  the case when the chain elements 
$\l_{q_i}$  of the tree-type diagrams 
have different number of edges $q_i\ge 1$. 


\section{Pr\"ufer code for tree-type diagrams}


The Pr\" ufer code has been developed to prove the Cayley formula for 
the number of free trees on $k$ vertices (see for instance \cite{Kaj}). There exists a wast  literature
concerning the use of the Pr\" ufer code and its versions and generalizations,
including various  types of trees \cite{BBR, Kaj}. 
In particular,  this method can be used  
to enumerate the trees with $k$ oriented labeled edges \cite{BBR}  with the help of the Cayley formula. 
What we do here resembles the reasoning of \cite{BBR}  with a version of the Pr\"ufer code
developed directly for the rooted trees with $k$ labeled edges (see \cite{BBR}, Lemma 2 and its proof).  
One could also trace out a relation with the paper \cite{Kaj}, where the 
Pr\" ufer code has been developed for the number of connected graphs assembled with $k$ blocks (graphs without cut points). However, no explicit formulas of the form (6) have been presented in these papers. 

\subsection{Coding   tree-type diagrams}




\noindent 
Regarding a diagram $D_k^{(q)}$,  it is convenient to identify the first and the last vertices of each element $\l_q$. In this cas we get 
a diagram constructed with the help of cycle elements 
$ \mu_q$ 
that have  one marked vertex and $q$ oriented simple edges. 
We denote such a tree-type diagram 
by $T_k^{(q)}$.
On Figure 2, we present a tree-type diagram $T_4^{(3)}$ obtained from 
the diagram $D_4^{(3)}$ depicted on Figure 1, where for convenience
we denoted $1$ by $a$, $2$ by $b$, etcetera.

\begin{figure}[htbp]
\centerline{\includegraphics
[height=10cm, width=16cm]{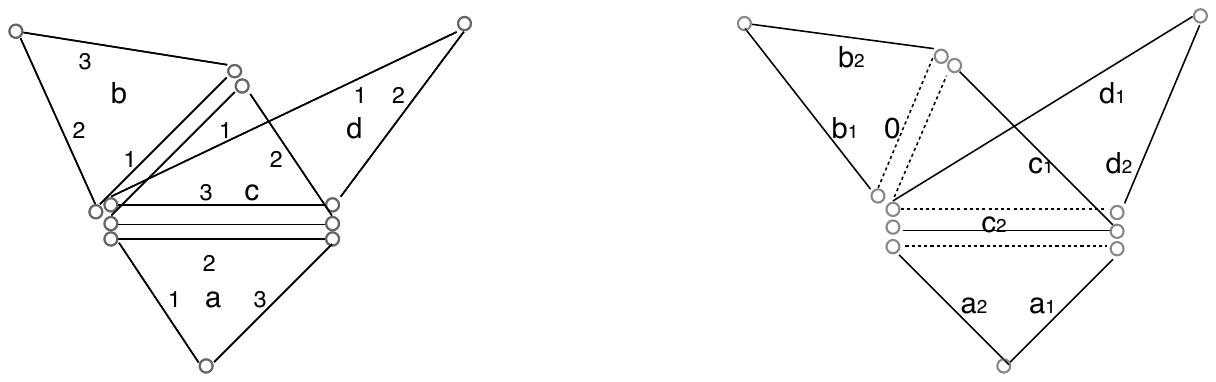}}
\caption{{ Tree-type diagrams $T_4^{(3)} $ and $\tilde T_4^{(3)}$}}
\end{figure}

\vs 
We assume  elements  $ \mu _1,  \mu _2, \dots   ,  \mu _k$ 
to be  colored in $k$ different colors; these colors   are alphabetically ordered, we denote them by 
$(a,b,c,\dots,h)$;
all edges of $  \mu_i$ are colored with the same color.

We are going to show that  each tree-type diagram $T_k^{(q)}$ 
generates 
 a sequence of $k-1$ symbols 
 $P_k^{(q)}= P(T_k^{(q)})$
that we refer to as the Pr\" ufer sequence of $T_k^{(q)}$. 
To write down this sequence, 
  we  modify $T_k^{(q)}$ 
as follows:

1) forget the order of edges in each element $  \mu_i$ of $T_k^{(q)}$;
 
2)  choose a couple of vertices $u,v$ of 
 $T_k^{(q)}$ that are connected by   simple or multiple edge that we denote by $(u,v)$; we say that this (multiple) edge is the root edge and attribute to it symbol $"0"$.
 Wash out  colors of  edges of elements
  $  \mu_{i_1}, \dots ,   \mu_{i_l}$ that join $u$ and $v$;
  all edges of  elements $  \mu_{i_1}, \dots ,  \mu_{i_l}$  
excepting the edges of the  root  $(u,v)$ keep their colors that they have had; we say that these color edges of  elements $  \mu_{i_1}, \dots ,   \mu_{i_l}$ are dominant edges;

3) consider an element $  \mu'_j$  attached to one of the dominant edges and wash out  the color of the edge of $  \mu'_j$ attached to this dominant edge;
all other edges of $  \mu'_j$ keep their color; these color edges are included to the list of the dominant edges;

4)  repeat the action (3) as many times as possible. 
\vs   
\noindent After this procedure performed, each element $  \mu_i$, $1\le i\le k$ has $q-1$ colored edges and one  colorless edge;  we refer to this colorless edge as to the bottom edge of 
$  \mu_i$. 
The diagram obtained $\hat T_k^{(q)}$ contains  $(q-1)k$ colored edges
and $k$ colorless edges.

\vs 5) regarding $\hat T_k^{(q)}$, consider
the left vertex of the bottom edge of each element $ \mu_j$  and  announce it to be the marked vertex; order the color edges of the same element
by labeling them with subscripts in the clockwise direction, say $a_1, a_2, \dots, a_{q-1}$. 
\vs 
\n Step (5) being performed, we get a diagram $\tilde T_k^{(q)}$. 
On Figure 2, we depict the diagram $\tilde T_4^{(3)}$ obtained from $T_4^{(3)}$.  Here
the colorless edges  are represented by dotted lines.

\vs 
With the diagram $\tilde T_k^{(q)}$, we are ready to write down 
a sequence  $P_k^{(q)}= P(\tilde T_k^{(q)})= P(T_k^{(q)})$ with the help of standard procedure
known for trees:

$\alpha$) find the minimal color such that there is no other element of
$\tilde T_{k}^{(q)}$  attached to the element $  \mu'$ of this color;

$\beta$) write down the color with subscript of the edge that $  \mu'$ is attached to; in the case when the element $  \mu'$ 
has a root edge, write down "0" instead of the color with subscript;

$\gamma$) erase the element $  \mu'$ from the diagram $T_k^{(q)}$ and 
denote the diagram obtained by $ \tilde T_{k-1}^{(q)}$;

$\delta$) perform the steps 
$(\alpha), (\beta)$ and $ (\gamma)$ to obtain  $\tilde T_{k-2}^{(q)}$, then $\tilde T_{k-3}^{(q)}$  and so on.
\vs 
\noindent We repeat  the action  $(\delta)$ till the last step  number $k-1$ is performed. On this stage 
the diagram $\tilde T_1^{(q)}$ contains only one element that has the root edge.
As a result,  we get a sequence of $k-1$ ordered cells (windows) fulfilled  by   letters with subscripts; 
 we denote this  sequence by $P_{k-1}^{(q)}$ and say that this is the Pr\"uder-type sequence (or code) of $T_k^{(q)}$.
Regarding the diagram $T_4^{(3)}$ depicted on Figure 2, we get the following sequence $P_3^{(3)}=P(T_4^{(3)})$:
\vs 
$$
P_3^{(3)}= \begin{tabular}{|c|c|c|c|}
 $c_2$& $ 0$ & $ c_2$
\\
\hline
\end{tabular} \ .
\eqno (7)
$$


\vs 
\vs

 \subsection{Construction of diagrams}

We determine a decoding procedure that is similar to  that used for the free trees:

-   start by writing down  a  Pr\" ufer-type sequence $P_{k-1}^{(q)}$
added by an auxiliary value $\{{\bf 0}\}$;   
under the sequence $P_{k-1}^{(q)}\sqcup \{{\bf 0}\}$ of $k$ symbols,    write down a line $L_k$ of $k$ symbols $a,b,c, \dots, h$;
we say that $(P_{k-1}^{(q)}\sqcup\{{\bf 0}\}, L_k)$ represents the  Pr\" ufer table;

i) in $L_k$,   underline  all letters  that appear in $P_{k-1}^{(q)}$;

ii)  
take the  value of the first window of $P_{k-1}^{(q)}$ that we denote by  $w_1$ 
and find the first value from the line $L_k$ that is not underlined, say $r$; 
then  join 
the element colored by  $r$  to the edge indicated by $w_1$; erasing $w_1$ and the first window  from $P_{k-1}^{(q)}$ and erasing 
$r$ from $L_k$ we get    the color Pr\" ufer table
$(P_{k-1}^{(q)}\sqcup\{{\bf 0}\}, L_{k-1})$.

\n Performing  the step  (ii) $k-1$ times, we arrive at the last stage when 
only one  element from the line $L_1$ remains. The last formal "step" is to join this element
with the root auxiliary edge marked by ${\bf 0}$.

\begin{figure}[htbp]
\centerline{\includegraphics
[height=10cm, width=16cm]{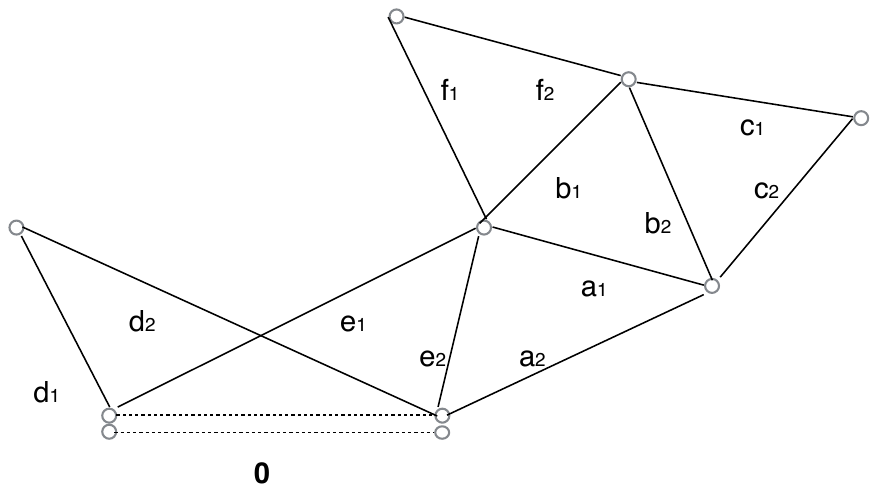}}
\caption{{ A tree-type diagram $\tilde T_6^{(3)} $}}
\end{figure}

For example, given a Pr\" ufer-type  sequence 
$P_5^{(3)}= (b_2, 0, b_1, a_1,e_2)$,
we get  corresponding table
$(P_5^{(3)} \sqcup \{{\bf 0}\}, L_6)$ as follows,
\[
 \begin{tabular}{c|c|c|c|c|c|c|c|c|c|c|c|c|} 
$ P_5^{(3)}\sqcup \{\bf 0\}$ & $b_2$ & $0$ & $b_1$ & $a_1$ & $e_2$ &  {\bf 0}
 \\ \hline
 $  L_6 $         &    \underline{a} & \underline{b} &c&  d & \underline{e}&  f        
 \end{tabular}\ .
 \eqno (8)
 \]
To construct $\tilde T_6^{(3)}$, we perform five principal steps and one formal step as follows:
$$
c\mapsto b_2, \quad d\mapsto 0, \quad f\mapsto b_1, \quad b\mapsto a_1, \quad a \mapsto e_2, \quad "e\mapsto {\bf 0}",
\eqno (9) 
$$
where the last term $"e\mapsto {\bf 0}"$ indicates the formal "step" 
saying that element $e$ is attached to  the root edge ${\bf 0}$.
Diagram $\tilde  T_6^{(3)}$ obtained form (8) is shown on Figure 3. For simplicity, we do not present here  dotted bottom edges, excepting the 
root edge marked by "0".

\vskip 0.2cm 

It is an easy exercice to show  that each $\tilde T^{(q)}_k$ generates a unique Pr\" ufer-type sequence $P_{k-1}^{(q)}$ 
and each Pr\"ufer sequence $P_{k-1}^{(q)}$  represents a unique tree-type diagram of the form $\tilde T_k^{(q)}$. 
There cardinality of the set of all possible  Pr\"ufer-type sequences
$P_k^{(q)}$ is given by 
$
|\CP_k^{(q)}|= (k(q-1)+1)^{k-1}
$
and therefore  
$$
| \tilde \CT_k^{(q)}| = (k(q-1)+1)^{k-1}.
\eqno (10)
$$

Eliminating the choice of the  root edge, we divide the number of rooted diagrams by $(q-1)k+1$. 
Rotating the special marked vertex of each element in $T_k^{(q)}$,
we get the factor $q^k$. 
This gives  (6). Theorem 1 is proved.
\hfill $\Box$  

\begin{figure}[htbp]
\centerline{\includegraphics
[height=8cm, width=16cm]{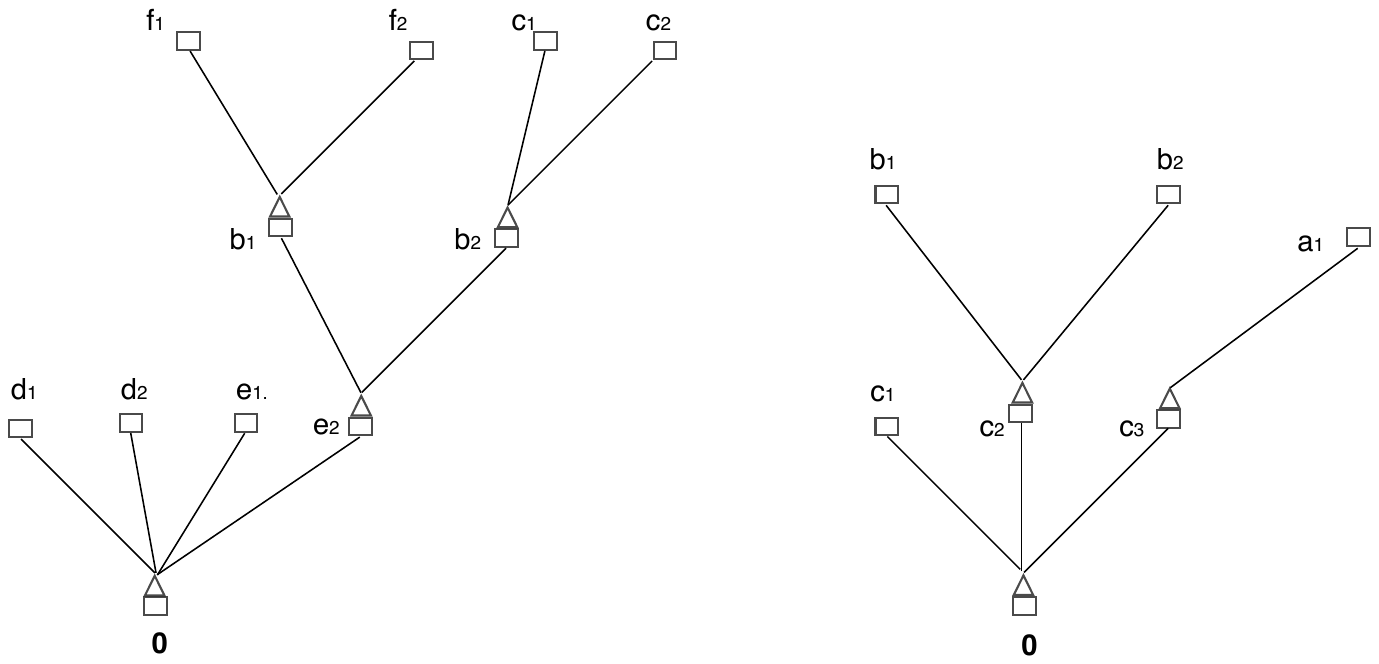}}
\caption{{ A regular tree diagram $\tilde \Gamma_6^{(3)} $ and irregular tree diagram $\tilde \Gamma^{(\{2,3,4\})}_3$}}
\end{figure}

Regarding the last graphical representation of $\tilde T_k^{(q)}$ 
given on Figure 3, one can observe that  
 that elements $  \mu_i$ of the diagram $T_k^{(q)}$ can also be represented as   a kind of stars;  
 each star of $\mu_i$  has  $q-1$ edges attached to a triangle vertex that corresponds 
 to the bottom edge of $\mu_i$. 
 we denote 
other  vertices of $q-1$ edges  by rectangles. 
 The rectangle  vertices can be used as a base for the triangle vertices to be put on. The auxiliary root vertex $"0$ is of rectangle form.
 The use of stars instead of elements $\mu_i$ 
means a passage to a kind of dual graph representation. 
 Finally, we get from $\tilde T_k^{(q)}$ a tree diagram $\tilde \Gamma_k^{(q)}$. 
 
 On Figure 4 we present $\tilde \Gamma^{(3)}_6$ that corresponds to 
 $\tilde T^{(3)}_6$ shown on  Figure 3. 
 Then relation between 
 tree-type diagrams $T_k^{(q)}$ and trees becomes even more clear.
 Let us note that in the case of $q=2$, diagrams  
 $\tilde \Gamma_k^{(2)}$ correspond to  ordinary rooted  
 trees  with $k$ labeled oriented edges. This gives a clear  
 explanation of  formula (5) obtained before by other authors, in particular   \cite{BBR}.

\section{Tree-type diagrams with edge multiplicity weight}

Let $V_i$, $i\in \bN$ be a given sequence of positif numbers. 
Given a tree-type diagram $T_k^{(q)}$, we  attribute to it a weight
$$
 W(T_k^{(q)}) = 
 \prod_{i=0}^{k-1} V_{i+1}^{s_i},
 \eqno (11)
$$
where $s_i$ is equal to the number of edges of multiplicity $i+1$. 
For example, the weight of  diagram $T_4^{(3)}$ shown on Figure 2 is 
equal to $W(T_4^{(3)})= V_1^7 V_2 V_3$.

\vs 
{\bf Theorem 2.} {\it The total sum of the set of weighted tree-type diagrams
$$
W_V(\CT_k^{(q)})= \sum_{T_k^{(q)}\in \CT_k^{(q)} }
W(T_k^{(q)})
$$
 is given by expression
$$
W_V(\CT_k^{(q)})= {q^k\over k(q-1)+1}
 \sum_{m=1}^{k-1}  \,  {{k(q-1)+1}\choose{m}}\ m! \ (k-1)! \
\
\sum_{ \stackrel{\s_{k-1}=(s_1, \dots, s_{k-1}), \ s_i\ge 0} 
{\| \s\|= k-1,  \, |\s| = m}}  \ 
\prod_{i=0}^{k-1} {V_{i+1}^{s_i}\over  ( i!)^{s_i}\, s_i!},
\eqno (12)
$$
where $\Vert \s_{k-1} \Vert = \sum_{i=1}^{k-1} i s_i
$, $| \s_{k-1}| = \sum_{i=1}^{k-1}  s_i$ and 
$s_0=kq- \|\s_{k-1}\|= k(q-1)+1$. }

\vs
{\bf Corollary.} Assuming that there exists a random variable $ \xi $ such that $V_{l+1}$ represents the $l$-th moment of $\xi$, 
$
V_{l+1} = \bE  \xi ^l 
$, $l\ge 0$,
we
deduce from (12) 
that
$$
W_V(\CT_k^{(q)})={q^k\over k(q-1)+1}
\bE   \left(  \xi_1+ \xi_2+ \dots+  \xi_{k(q-1)+1}\right)^{k-1},
\eqno (13)
$$
where $( \xi_1,  \xi_2, \dots ,  \xi_{k(q-1)+1})$ is the family of 
jointly independent random variables that have the same probability distribution as $\xi $ 
(see \cite{K-24} for more details).
\vs
Proof of Theorem 2 is based on a simple observation that if 
the Pr\"ufer-type sequence $P_{k-1}^{(q)}$ contains exactly $i$ cells with the same symbol, then corresponding diagrams $\tilde T_k^{(q)}$ and 
$T_k^{(s)}$ contain a multiple edge of multiplicity $i+1$ (see example 
$P_3^{(3)}$ (7) and corresponding diagram $T_4^{(3)}$ on Figure 2).

\vs 
Regarding variable   $\s_{k-1}=(s_1,s_2,\dots, s_{k-1})$, $0  \le  s_i  \le k-1$,
we say  that a  sequence $P_{k-1}^{(q)}$ is of $\s_{k-1}$-type,
with 
  if it contains 
  $s_1$ subsets of one element, $s_2$ subsets  of two elements, 
  $s_{k-1}$ subsets of $k-1$ elements, in each subset the  elements 
  are identical. 
  It is clear that the number of different values that 
  appear in $\s$-type sequence $\CP_k^{(q)}$ is given by 
  $| \s |= s_1+s_2+\dots+ s_k$.

Regarding $\s_{k-1}$ such that $|\s_{k-1}|= m$, we have to choose
from $k(q-1)+1$ variables $m$ values. This can be done by 
${{k(q-1)+1}\choose{m}}m!$ ways.    
  Then we can write the following identity for the number of all possible Pr\" ufer  sequences 
  
$$
 \sum_{m=1}^{k-1}  \,  {{k(q-1)+1}\choose{m}}\ m! \ 
\ \sum_{ \stackrel{\s_{k-1}=(s_1, \dots, s_{k-1}), \ s_i\ge 0} 
{\| \s_{k-1}\|= k-1,  \, |\s_{k-1}| = m}} \ (k-1)! \
\prod_{i=1}^{k-1} {1\over  ( i!)^{s_i}\, s_i!}= 
\big( k(q-1)+1\big)^{k-1}.
\eqno (14)
$$ 
Then relation (12) is obtained as an obvious generalization of  (13).
Theorem 2 is proved. \hfill $\Box$
\vs

\section{Tree-type diagrams with non-regular chains}

 In this section we consider non-regular tree-type diagrams 
 $D^{(\{q_1, \dots, q_k\})}_k$ assembled from  oriented chain  
  elements $  \l_1, \dots,   \l_k$ that have $q_1, \dots, q_k$ edges, respectively.
  On Figure 4 we present an exemple of an irregular tree diagram $\tilde \Gamma^{(\{2,3,4\})}_3$.
 \vs 
{\bf Theorem 3.} {\it Given $k\in \bN $ and a sequence of naturals $\{q_1, q_2, \dots, q_k\}$, $q_i\ge 2$, the cardinality  of the family of all tree-type diagrams 
assembled from oriented  chain elements $\l_{q_1}, \l_{q_2}, \dots, \l_{q_k}$ is given by 
 $$
| \CD_k^{(\{q_1, q_2, \dots, q_k\})}|= 
d^{(q_1,q_2, \dots, q_k)}_k=  \big( q_1+q_2+ \dots + q_k - k+ 1)^{k-2} \ \prod_{i=1}^k q_i\ . 
 \eqno (15)
 $$
Remark. } Relation (15) remains true in the case when 
$q_j=1$ for one or several $j$, $1\le j\le k$. 
 \vs 
 The proof of Theorem 3 repeats the proof of Theorem 1 with the only difference that the set of the Pr\"ufer-type sequences of length $k-1$
 is given by $(q_1+ q_2+ \dots q_k - k+1)^{k-1}$. 
 

Regarding irregular tree-type diagrams $D_k^{(\{q_1,q_2, \dots, q_k\})}$ with the  weight (11),
$$
W\big(D_k^{(\{q_1,q_2, \dots, q_k\})}\big)= \prod_{i=0}^{k-1} V_{i+1}^{s_i},
 $$
 we observe that similar to (12) equality
 is true, where $q^k$ is replaced by $\prod_{i=1}^k q_i$,
 factor $k(q-1)+1$ is replaced by 
 $(q_1+q_2+ \dots + q_k - k+1)$ and the binomial coefficient
 ${{k(q-1)+1}\choose{m}}$ is replaced by the binomial coefficient 
 ${{q_1+q_2+\dots + q_k - k+1}\choose{m}}$. As a consequence,
 we obtain analog of relation (13),
 $$
W_V\big(\CD_k^{(\{q_1,q_2, \dots, q_k\})}\big)={q_1 q_2 \cdots q_k
\over q_1+ q_2 + \dots + q_k-k+1}
\bE   \left(  \xi_1+ \xi_2+ \dots+  \xi_{q_1+ q_2+ \dots + q_k - k+1}\right)^{k-1}.
\eqno (16)
$$

Relation (16), as well as (13) can be useful in the studies of the
asymptotic properties of  
weighted diagrams in the limiting transition of 
infinite $q$ or $q_i$ (see \cite{K-24} for more details). We do not comment on this topic here.

\vs 
\section{Appendix. P\'olya equation and Lagrange formula}


In paper \cite{K-08}, it is shown that the  generating function, 
$$
H_q(z) = \sum_{k\ge 0}h_k^{(q)}  z^k,
\quad h_k^{(q)}= (q-1)k+1) d_k^{(q)}/k!
$$
verifies equality
$$
H_q(z) =
\exp  \left\{ qz (H_q(z))^{q-1}\right\}
\eqno (17)
$$
known  as the P\'olya equation \cite{Sta}. 
This relation  has been deduced from  a recurrence that determines numbers 
$d_k^{(q)}$ obtained in \cite{K-08} (see also \cite{FU}).  
However, equation (17) has been solved in \cite{K-08} in the case of $q=2$ only.
Let us describe briefly its solution in the general case of 
$q\ge 2$ 
obtained with the help of the contour integration method 
of the Lagrange inversion theorem \cite{PS}. 

Let us consider an auxiliary  function 
$$
t(w)= w e^{-q(q-1)w}, \quad w \in \bC
\eqno (18)
$$
that is analytic in a vicinity of zero with the  derivative  non-zero in the domain
$\{w\in \bC: |w |< (q(q-1))^{-1}\}$ 
and has its inverse $t^{-1}(z)= \psi(z)$. 
Then $\psi(z)$, $z \in \bC$  is also analytic in a vicinity of zero. The same is true for 
$$
H_q(z)= e^{q \psi(z)} = \sum_{k=0}^\infty C_k z^k,
$$
where the series converges and where $C_k= h_k^{(q)}$. 
Then  we can write that 
$$
C_k = {1\over 2\pi i} \oint {e^{q\psi(z)}\over z^{k+1}} \, dz.
\eqno (19)
$$
Taking into account (18), relation $z= t(w)$ and equality 
$$
dz =   \left( e^{-q(q-1)w} - q(q-1) w e^{-q(q-1)w} \right)\, dw
$$
that follows from (18), we rewrite  (19) in the form
$$
C_k = {1\over 2\pi i} \oint { e^{qw + q(q-1)wk} \over w^{k+1} }  
\left(
1 - q(q-1)w\right) dw
$$
$$
=   \lfloor e^{q(q-1)wk + qw}\rfloor_{(k)}
-
q(q-1)    \lfloor e^{q(q-1)wk + qw}\rfloor_{(k-1)},
$$
where we denoted by $  \lfloor f(w)\rfloor_{(k)}$ the $k$-th coefficient of the Taylor expansion of $f(w)$ at zero. 
Then 
$$
h^{(q)}_k = {q^k\over k!} \big( (q-1)k+1\big)^{k-1} - 
{q^k(q-1)\over (k-1)!} \big( (q-1)k+1\big)^{k-1}= 
{q^k\over k!} \big( (q-1)k+1\big)^{k-1}.
$$
Remembering definition of $h^{(q)}_j$, we get the final expression for $d_k^{(q)}$ (6). This result confirms Theorem 1. However, the use of 
generating function technique can be much more difficult in the case of weighted and irregular tree-type diagrams, where the Pr\"ufer codification
is fairly simple and efficient.

\end{document}